\documentclass{amsart}
\usepackage{setspace, amssymb, amsmath, amsfonts}

\theoremstyle{plain}
\newtheorem{thm}{Theorem}[section]
\newtheorem{theorem}[thm]{Theorem}

\newtheorem{lemma}[thm]{Lemma}

\theoremstyle{definition}

\theoremstyle{remark}

\newcommand{\ld}{\lambda}

\newcommand{\R}{\mathbb{R}}
\newcommand{\Z}{\mathbb{Z}}

\DeclareMathOperator{\Aut}{Aut}
\DeclareMathOperator{\aia}{AIA}

\newcommand{\semi}{\ltimes}

\begin{document}

\title{An isospectral deformation on an orbifold quotient of a nilmanifold}

\author[E. Proctor]{Emily Proctor}
\address{Emily Proctor \\ Middlebury College, Department of Mathematics, Middlebury, VT, 05753}
\email{eproctor@middlebury.edu}
\author[E. Stanhope]{Elizabeth Stanhope}
\address{Elizabeth Stanhope \\ Lewis and Clark College, Department of
Mathematical Sciences, 0615 SW Palatine Hill Road, MSC 110, Portland, OR 97219}
\email{stanhope@lclark.edu}
\maketitle

\begin{abstract}
We construct a Laplace isospectral deformation of metrics on an orbifold quotient of a nilmanifold.  Each orbifold in the deformation contains singular points with order two isotropy.  Isospectrality is obtained by modifying a generalization of Sunada's Theorem due to DeTurck and Gordon.
\end{abstract}

\medskip

\noindent 
\begin{center}
\begin{small}
\end{small}
\end{center}
\bigskip

\section{Introduction}

A Riemannian orbifold (see \cite{S}, \cite{Sc}) is a mildly singular generalization of a Riemannian manifold.  For example the quotient space of a Riemannian manifold under an isometric, properly discontinuous group action is a Riemannian orbifold \cite{T}.  First defined in 1956 by I. Satake, orbifolds have proven useful in many settings including the theory of 3-manifolds, symplectic geometry, and string theory.  

The local structure of a Riemannian orbifold is given by the orbit space of a Riemannian manifold under the isometric action of a finite group.  If a point, $p$, in the manifold is fixed under the group action, the corresponding element of the orbit space, $\bar{p}$, is called a \emph{singular point} of the orbifold.  The \emph{isotropy type} of a point $\bar{p}$ in the orbit space is the isomorphism class of the isotropy group of a point $p$ in the manifold that projects to $\bar{p}$ under the quotient.  The \emph{singular set} of an orbifold is the set of all  singular points of the orbifold.

The tools of spectral geometry can be transferred to the setting of Riemannian orbifolds by exploiting the well-behaved local structure of these spaces (see \cite{C}, \cite{St}).  As in the manifold setting, the spectrum of the Laplace operator of a compact Riemannian orbifold is a sequence
$
0 \le \ld_0 \le \ld_1 \le \ld_2 \le\dots \uparrow +\infty
$ 
where each eigenvalue has finite multiplicity.  We say that two orbifolds are \emph{isospectral} if their Laplace spectra agree.

In this note we show that the formulation of Sunada's Theorem found in \cite{DG} can be used to obtain isospectral deformations on Riemannian orbifolds with non-trivial singular sets.  We prove this fact in Section 2 by observing that the proof of Theorem 2.7 in \cite{DG} does not require that the action of the discrete subgroup $\Gamma$ be free.  In Section 3 we display an example of an isospectral deformation of metrics on an orbifold quotient of a nilmanifold.    

The only other known examples of non-manifold isospectral deformations on orbifolds were recently obtained by Sutton using a blend of the torus action method and the Sunada technique \cite{Su}.  Other examples of non-manifold isospectral orbifolds include pairs with boundary in \cite{BCDS}, and in \cite{BW}; isospectral flat 2-orbifolds that are not conjugate (in terms of lengths of closed geodesics) \cite{DR}; a $(2m)$-manifold isospectral to a $(2m)$-orbifold on $m$-forms \cite{GR}; pairs of isospectral orbifolds for which the maximal isotropy groups have different orders \cite{R} and \cite{Sch}; and arbitrarily large finite families of isospectral orbifolds in \cite{SSW}. 

\subsection{Acknowledgments} The authors thank Carolyn S. Gordon for her helpful suggestions during the course of this work.

\section{Isospectral deformations on orbifolds}

\bigskip

In this section we observe that the generalization of Sunada's method found in \cite{DG} can be further generalized to include isospectral deformations of metrics on orbifolds.  In what follows we will assume that $G$ is a Lie group with simply connected identity component $G_0$. We let $\Gamma$ be a discrete subgroup of $G$ such that $G=\Gamma G_0$ and $(G_0\cap \Gamma)\backslash G_0$ is compact.  

Given an automorphism $\Phi: G\to G$, we say that $\Phi$ is an \textit{almost-inner automorphism} if for each $x\in G$ there exists an element $a\in G$ such that $\Phi(x) = axa^{-1}$.  More generally, if $\Phi: G\to G$ is an automorphism such that for each $\gamma\in\Gamma$ there exists $a\in G$ satisfying $\Phi(\gamma) = a\gamma a^{-1}$, we say that $\Phi$ is an \textit{almost-inner automorphism of $G$ relative to $\Gamma$}.  We denote the set of all almost-inner automorphisms of $G$ (resp. almost-inner automorphisms of $G$ with respect to $\Gamma$) by $\mathrm{AIA}(G)$ (resp. $\mathrm{AIA}(G;\Gamma)$).

We have the following theorem.

\begin{theorem} \label{DG}\cite{DG} Let $G$, $G_0$, and $\Gamma$ be as above with $G_0$ nilpotent and let $\Phi \in \mathrm{AIA}(G;\Gamma)$.  Suppose that $G$ acts effectively and properly on the left by isometries on a Riemannian manifold $(M,g)$ and that $\Gamma$ acts freely on $M$ with $\Gamma\backslash M$ compact.  Then, letting $g$ denote the submersion metric, $(\Phi(\Gamma)\backslash M, g)$ is isospectral to $(\Gamma\backslash M, g)$.  
\end{theorem}

The proof of Theorem~\ref{DG} is based on work by Donnelly in \cite{D} concerning the existence of a heat kernel on a manifold $M$ which admits a properly discontinuous (but not necessarily free) action by a group $\Gamma$.  Donnelly shows that if $\Gamma\backslash M$ is compact, then there exists a unique heat kernel on $M$.  Furthermore, Donnelly gives the following relationship between the heat kernels on $M$ and on $\Gamma\backslash M$. 

\begin{theorem}\label{D}\textrm{\cite{D}} Let $\Gamma$ act properly discontinuously on $M$ with compact quotient $\overline{M} = \Gamma\backslash M$.  Suppose that $F$ is a fundamental domain for $\Gamma\backslash M$.  If $\bar{x}, \bar{y}\in\overline{M}$ then set
\begin{equation*}
\overline{E}(t,\bar{x},\bar{y}) =\sum_{\gamma\in\Gamma} E(t,x,\gamma\cdot y)
\end{equation*}
where $x,y\in F$, $\bar{x}=\pi(x)$, and $\bar{y}=\pi(y)$.  If $E$ is the heat kernel of $M$, the sum on the right converges uniformly on $[t_1,t_2]\times F \times F$, $0<t_1\leq t_2$, to the heat kernel on $\overline{M}$.
\end{theorem}

Notice that since the action of $\Gamma$ need not be free, the quotient space $\overline{M}$ may not be a manifold.

Theorem~\ref{DG} relies on the fact that two manifolds $(M_1,g_1)$ and $(M_2,g_2)$ are isospectral if and only if they have the same heat trace, i.e. $\int_{M_1}E_1(t,x,x)\,dx = \int_{M_2}E_2(t,x,x)\,dx$, where $E_i$ denotes the heat kernel on $M_i$.  In particular, the proof uses Theorem~\ref{D} to pull the heat trace back from the quotient $\Gamma\backslash M$ to the cover $M$ in order to use combinatorial arguments to reexpress the heat trace on $\Gamma\backslash M$.  The new expression of the heat trace makes it evident that, when comparing the heat trace of $(\Gamma\backslash M,g)$ with the heat trace of $(\Phi(\Gamma)\backslash M, g)$, if certain volumes (which depend only on $\Gamma$ and $\Phi(\Gamma)$) are equal then the respective heat traces are equal.  DeTurck and Gordon show that when $\Phi$ is an almost-inner automorphism, these volumes are in fact equal, and hence $(\Gamma\backslash M, g)$ and $(\Phi(\Gamma)\backslash M, g)$ are isospectral. 

We note that, as with Theorem~\ref{D}, the proof of Theorem~\ref{DG} does not rely on the freeness of the action of $\Gamma$ on $M$.  Therefore we make the following generalization of Sunada's theorem.

\begin{theorem}\label{ourtheorem}
Suppose that $G$, $G_0$, and $\Gamma$ are as above and $G_0$ is nilpotent.  Suppose that $G$ acts effectively and properly discontinuously on the left by isometries on $(M,g)$ with $\Gamma\backslash M$ compact.  Let $\Phi\in\aia(G;\Gamma)$.  Then, letting $g$ denote the submersion metric, the quotient orbifolds $(\Gamma\backslash M, g)$ and $(\Phi(\Gamma)\backslash M, g)$ are isospectral. 
\end{theorem}

\section{Examples}

Now we apply Theorem~\ref{ourtheorem} to give an example of a nontrivial isospectral deformation on an orbifold. We first note the following.

\begin{lemma}  \label{isometry} 
If $\Phi$ is an isomorphism of $G$ and $\Gamma$ acts on $M$ in such a way that there exists a diffeomorphism $\Psi$ of $M$ satisfying $\Psi(g\cdot x) = \Phi(g)\cdot\Psi(x)$ for all $g\in G$ and $x\in M$, then $(\Phi(\Gamma)\backslash M, g)$ is isometric to $(\Gamma\backslash M, \Psi^* g)$.  
\end{lemma}

\noindent Applying Theorem~\ref{ourtheorem} in conjunction with Lemma \ref{isometry} will allow us to produce an isospectral deformation on a fixed orbifold $\Gamma\backslash M$.

In Appendix B to \cite{DG}, K. B. Lee translates Theorem~\ref{DG} to the setting of infranilmanifolds.  For a group $G$ we have that $\Aut(G)\semi G$ acts on $G$ by $(\phi,g)\cdot h = g\phi(h)$.  Consider the case when $G$ is a simply connected nilpotent Lie group and $\Gamma$ is a uniform discrete subgroup of $G$.  Take $\Pi$ to be a finite extension of $\Gamma$ in $\Aut(G)\semi G$.  If the action of $\Pi$ on $G$ is free, then $\Pi\backslash G$ is an infranilmanifold.    Lee observes that by setting $\Gamma$, $G_0$, and $G$ from Theorem~\ref{DG} equal to $\Pi$, $G$ and $\Pi G$, and assuming that the action of $\Pi$ on $G$ is free, we can find isospectral deformations on infranilmanifolds.   We note that a priori, the action of $\Pi G$ on $G$ need not be free.  Thus by working in this setting we introduce the possibility of finding isospectral orbifold quotients of $G$.

Lee gives a specific example to illustrate his case.  His example is based on a similar example found in \cite{GW2}.  

Let $G$ be the Lie group 
\begin{equation*}
\{(x_1, x_2, y_1, y_2, z_1, z_2)\vert x_i, y_i, z_i\in \mathbb{R}\}
\end{equation*}
where group multiplication is defined by 
\begin{multline*}
(x_1,\dots,z_2)(x_1^{\prime},\dots z_1^{\prime}) \\
=(x_1 + x_1^{\prime}, \dots, y_2+y_2^{\prime}, z_1+z_1^{\prime} +x_1y_1^{\prime} +x_2y_2^{\prime}, z_2 +z_2^{\prime} + x_1y_2^{\prime}).
\end{multline*}
Suppose that $\Gamma$ is the integer lattice in $G$ and define $\Phi_t: G\to G$ by 
\begin{equation*}
\Phi_t(x_1,x_2,y_1,y_2,z_1,z_2) = (x_1,x_2,y_1,y_2,z_1,z_2+ty_2),
\end{equation*}
where $t\in [0,1)$.  In the original example Gordon and Wilson show that each $\Phi_t$ is an almost-inner automorphism so, applying Lemma \ref{isometry},  $(\Phi_t(\Gamma)\backslash G, g)$ is an isospectral deformation on $\Gamma\backslash G$.  They also show that the deformation is nontrivial.

In his example, Lee defines $\alpha\in \Aut(G)\semi G$ by 
\begin{equation*}
\alpha(x_1,x_2,y_1,y_2,z_1,z_2) = (x_1, x_2, -y_1,-y_2,-z_1,-z_2+\tfrac{1}{2})
\end{equation*}
and lets $\Pi=\Gamma\cup\alpha\Gamma$.  Since $\alpha$ commutes with $\Phi_t$ for all $t$, we can 
extend each $\Phi_t$ to an element $\tilde\Phi_t$ of $\aia(\Pi G;\Pi)$.  If $g$ is a $\Pi G$-invariant metric on $G$, then for each $t$, $(\tilde\Phi_t(\Pi)\backslash G, g)$ is isospectral to $(\Pi \backslash G,g)$.  

Lee implicitly assumed that the action of $\Pi$ on $G$ is free.  However, we can see by closer inspection that the action of $\Pi$ on $G$ is not free.  For example, any point of the form $(x_1, x_2, 0,0,0,\frac{1}{4})$ is fixed by $\alpha\in\Pi G$.  In fact the set of all fixed points of the action of $\Pi$ on $G$ is:
\[ \{(x_1, x_2, y_1, y_2, z_1, z_2) \in \R^6 : x_1, x_2 \in \R, y_1, y_2, z_1 \in \tfrac{1}{2}\Z, z_2 = \tfrac{n}{2} + \tfrac{1}{4}\}\]
where $n$ is any integer.  The isotropy group of a point in this set has the form, 
\[\{1, (\Phi, (0, 0, 2y_1, 2y_2, 2z_1, 2z_2))\}.\]
So we see that $\Pi\backslash G$ is an orbifold containing singular points with $\Z_2$ isotropy type.  Thus Lee's example is an illustration of Theorem~\ref{ourtheorem}: after applying Lemma \ref{isometry} with $\Psi = {\Phi}_t$ and $\Phi = \tilde\Phi_t$, we have an isospectral deformation of metrics on the \textit{orbifold} $\Pi\backslash G$.

This example is a nontrivial deformation.  Indeed, suppose that $\tau: (\Pi\backslash G, g)\to(\Pi\backslash G, \Phi_t^*g)$ is an isometry.  Then because $G$ is simply connected and $\Pi$ is discrete, $G$ is the simply connected cover of $\Pi\backslash G$.  Thus $\tau$ lifts to an isometry, also called $\tau$ from $(G, g)$ to $(G, \Phi_t^*g)$.  Since $G$ is a nilpotent Lie group $\tau$ must be an element of $\Aut(G)\semi G$ (see \cite{GW1}).  Furthermore, because $\tau$ is a lift we have that $\tau\circ\Pi\circ\tau^{-1} = \Pi$ within the transformation group $\Aut(G)\semi G$.  On the other hand, $G$ is normal in $\Aut(G)\semi G$ so conjugation by $\tau$ maps $G$ to itself.  Therefore, conjugation by $\tau$ leaves $\Gamma$ invariant.  This implies that $\tau$ must descend to an isometry $\tau: (\Gamma\backslash G, g) \to (\Gamma\backslash G, \Phi_t^*g)$.  However, from \cite{GW2} we know that no such isometry can exist.  Thus $(\Pi\backslash G, g)$ cannot be isometric to $(\Pi\backslash G, \Phi_t^*g)$.  

Note that Lee's example can be modified to produce an isospectral deformation on a manifold.  
For example, suppose that we define $\beta: G\to G$ by 
\begin{equation*}
\beta(x_1, x_2,y_1,y_2,z_1,z_2) = (x_1,x_2,y_1,y_2, -z_1, z_2+\tfrac{1}{2}).
\end{equation*}
Letting $\Pi = \Gamma\cup\beta\Gamma$ we see that since $\beta^2$ is simply translation by $(0,0,0,0,0,1)$, $\Pi$ is a finite extension of $\Gamma$. Since $\beta$ commutes with the maps $\Phi_t$ defined above, we can extend each $\Phi_t$ to an element $\tilde\Phi_t$ of $\aia(\Pi G; \Pi)$.  Finally by direct computation we can see that the action of $\Pi$ on $G$ has no fixed points.  Thus we have an isospectral deformation of metrics on a manifold.  The proof that the deformation is nontrivial is identical to the one given above.


\begin{thebibliography}{BCDS}

\bibitem[BCDS]{BCDS} Buser, P., Conway, J., Doyle, P. and Semmler, K., \emph{Some planar isospectral domains}, Internat. Math. Res. Notices. \textbf{9} (1994), 391ff., approx.\ 9 pp.\ (electronic).

\bibitem[BW]{BW} B{\'e}rard, P. and Webb, D., \emph{On ne peut pas entendre l'orientabilit\'e d'une surface}, C. R. Acad. Sci. Paris S\'er. I Math. \textbf{320} (1995), no. 5, 533--536.

\bibitem[C]{C} Chiang, Y.-J., \emph{Spectral geometry of $V$-manifolds and
its application to harmonic maps}, in ``Differential geometry: partial
differential equations on manifolds" (Los Angeles, CA, 1990),  93--99,
Proc. Sympos. Pure Math., \textbf{54}, Part 1,
Amer. Math. Soc., Providence, RI, 1993.

\bibitem[D]{D} Donnelly, H., \emph{Asymptotic expansions for the compact quotients of properly discontinuous group actions}, Illinois J. Math. \textbf{23} (1979) 485--496.

\bibitem[DG]{DG} DeTurck, D. and Gordon, C., \emph{Isospectral deformations II}, Comm. Pure and Appl. Math. \textbf{42} (1989) 1067--1085 (with an appendix by K. B. Lee).

\bibitem[DR]{DR} Doyle, P. and Rossetti, J., \emph{Isospectral hyperbolic surfaces having matching geodesics}, preprint, ArXiv math.DG/0605765.

\bibitem[GR]{GR} Gordon, C. S., and Rossetti, J., \emph{Boundary volume and length spectra of {R}iemannian manifolds:  what the middle degree {H}odge spectrum doesn't reveal}, Ann. Inst. Fourier, \textbf{53} (2003), no. 7, 2297--2314.

\bibitem[GW1]{GW1} Gordon, C. and Wilson, E., \emph{Isometry groups of Riemannian solvmanifolds}, Trans. AMS, \textbf{307} (1988) 245--269.

\bibitem[GW2]{GW2} Gordon, C. and Wilson, E., \emph{Isospectral deformations of compact solvmanifolds}, J. Diff. Geom. \textbf{19} (1984) 241--256.

\bibitem[R]{R} Rossetti, J., private communication.

\bibitem[S]{S} Satake, I., \emph{On a generalization of the notion of
manifold}, Proc. Nat. Acad. Sci. U.S.A. \textbf{42} (1956), 359--363.

\bibitem[Sch]{Sch} Schueth, D., private communication.

\bibitem[Sc]{Sc} Scott, P., \emph{The geometries of $3$-manifolds}, Bull. London
Math. Soc. \textbf{15} (1983), no. 5, 401--487.

\bibitem[SSW]{SSW} Shams, N., Stanhope, E., and Webb, D., \emph{One cannot hear orbifold isotropy type.}, Arch. Math., \textbf{87} (2006), no. 4, 375--385.

\bibitem[St]{St} Stanhope, E. \emph{Spectral bounds on orbifold isotropy},
Annals of Global Analysis and Geometry \textbf{27} (2005), no. 4, 355--375.

\bibitem[Su]{Su} Sutton, C. J., \emph{Equivariant isospectrality and isospectral deformations of metrics on spherical orbifolds}, preprint, ArXiv  math.DG/0608557.

\bibitem[T]{T} Thurston, W., \emph{The Geometry and Topology of
Three-Manifolds}, lecture notes, Princeton University, Princeton, NJ, 1979.

\end{thebibliography}
\end{document}